\documentclass[12pt]{article}
\usepackage{amsfonts,amssymb,epsfig}
\newcommand{\ol}{\setlength{\itemsep}{0pt.}\begin{enumerate}}
\newcommand{\eol}{\end{enumerate}\setlength{\itemsep}{-\parsep}}
\date{}
\begin{document}
\newtheorem{df}{Definition}
\newtheorem{thm}{Theorem}
\newtheorem{lm}{Lemma}
\newtheorem{pr}{Proposition}
\newtheorem{co}{Corollary}
\newtheorem{re}{Remark}
\newtheorem{note}{Note}
\def\R{{\mathbb R}}

\def\E{\mathbb{E}}
\def\calF{{\cal F}}
\def\N{\mathbb{N}}
\def\calN{{\cal N}}
\def\calH{{\cal H}}
\def\n{\nu}
\def\a{\alpha}
\def\d{\delta}
\def\t{\theta}
\def\l{\lambda}
\def\e{\varepsilon}
\def\pf{ \noindent {\bf Proof: \  }}
\def\trace{\rm trace}\newcommand{\qed}{\hfill\vrule height6pt  width6pt
depth0pt}
\def\endpf{\qed \medskip}
\def\colon{{:}\;}
\setcounter{footnote}{0}

\title{Special orthogonal splittings of $L_1^{2k} $}

\author{Gideon Schechtman\thanks{Supported by the
Israel Science Foundation }} \maketitle

\begin{abstract}
We show that for each positive integer $k$ there is a $k\times k$
matrix $B$ with $\pm 1$ entries such that putting $E$ to be the
span of the rows of the $k\times 2k$ matrix $[\sqrt{k}I_k,B]$,
then $E,E^{\bot}$ is a Kashin splitting: The $L_1^{2k}$ and the
$L_2^{2k}$ are universally equivalent on both $E$ and $E^{\bot}$.
Moreover, the probability that a random $\pm 1$ matrix satisfies
the above is exponentially close to $1$.

\end{abstract}

\section{Introduction}
For $0<p<\infty$ and $n\in\N$ let $L_p^n$ denote $\R^n$ with the
norm $\|x\|_{L_p^n}=\Big(n^{-1}\sum_{i=1}^n |x_i|^p\Big)^{1/p}$,
where $x=(x_1,x_2,\dots, x_n)$. A celebrated theorem of Kashin
\cite{ka} states that $L_1^{2k}$ can be decomposed into two
orthogonal (with respect to the inner product induced by
$\|\cdot\|_{L_2^{2k}}$) $k$-dimensional subspaces on each of
which the two norms $\|\cdot\|_{L_1^n}$ and $\|\cdot\|_{L_2^n}$
are universally equivalent, i.e., putting, for a subset
$E\subseteq L_1^n$ \
\[
C_n(E)=\sup\{\|x\|_{L_2^n}/\|x\|_{L_1^n}\ ;\ x\in E, \ x\not=0 \},
\]
we can find a $k$-dimensional subspace of $L_1^{2k}$ for which
\[
C_{2k}(E),C_{2k}(E^{\bot})\le C
\]
where $C<\infty$ is some universal constant. We shall call such a
choice of (orthogonal) subspace(s) a Kashin splitting with
constant $C$.

The proof(s) of Kashin theorem are probabilistic and do not
produce an explicit subspace $E$ as above. For example it is shown
that with high probability over the orthogonal group $O(k)$ (where
the probability is the Haar measure), the span of the rows of the
$k\times 2k$ matrix $[U_1,U_2]$, where $U_1,U_2\in O(k)$ are
chosen independently, is such a subspace.

In a recent paper Anderson \cite{a} found an explicit
determinantal formula for $C_{2k}(E)$ and $C_{2k}(E^{\bot})$ for
any $k$-dimensional subspace $E$ (involving determinants of
$k\times k$ submatrices of the $k\times 2k$ matrix whose rows are
any basis of $E$). Anderson then continues and presents a
discretization of (a variant of) the random decomposition,
reducing the search of a Kashin splitting to a search among
$k\times 2k$ matrices with integer entries (ranging in some
bounded, though of size exponential in $k$, set). The point is
that this suggests a possibility of finding an {\em explicit}
Kashin splitting. It also permits a (not very efficient) search
algorithm for finding a good Kashin splitting (although, it
seems, for that the main point in the paper, the determinantal
formulas, can be avoided).

In this paper we take this direction one step farther by showing
that one can replace the integral matrices by matrices whose
entries are taken only from the set $\{0,\sqrt{k},1,-1\}$. More
precisely, we show in Theorem \ref{thm:main} that with high
probability for a random choice of $k\times k$ matrix $B$ with
entries being independent Bernoulli $\pm 1$ variables, the span of
the rows of the matrix $[\sqrt{k}I,B]$ form a Kashin splitting
with some universal constant. (Since $\sqrt{k}$ is not
necessarily an integer, one may wonder whether this is, strictly
speaking, a strengthening of Anderson's result. However, in the
proof of Theorem \ref{thm:main} bellow one can easily replace
$\sqrt{k}$ with $[\sqrt{k}]$ everywhere and get such a formal
strengthening.)

We would like next to indicate what Anderson's determinantal
formula gives for such matrices. For $R,C$ two subsets of
$\{1,2,\dots,k\}$ denote by $B_{R,C}$ the sumatrix of $B$ formed
by the rows in $R$ and the columns in $C$. For a submatrix
$D=B_{R,C}$ of $B$ and row $i\in R$ let $D_{-i}$ be the matrix
$B_{R\setminus\{i\},C}$ similarly for a $j\notin C$ let $D^{+j}$
be the matrix $B_{R,C\cup\{j\}}$. For $l=1,2,\dots,k-1$, a
$(l+1)\times l$ submatrix $D=B_{R,C}$ of $B$ and $p=1,2$ denote
\[
\Delta_p(D,B)=\Big(k^{p/2}\sum_{i\in R}|\det
D_{-i}|^p+\sum_{j\notin C}|\det D^{+j}|^p\Big)^{1/p}.
\]
Using Anderson's determinantal formulas one gets, as we shall
indicate in Corollary \ref{co:determinants}, that for a $k\times
k$ matrix $B$ the rows of $[\sqrt{k}I,B]$ form a Kashin
decomposition with constant
\[
\sqrt{2k}\max\Big\{\max\Big\{\frac{\Delta_2(D,B)}{\Delta_1(D,B)}\Big\},\max\Big\{\frac{\Delta_2(D,B^*)}{\Delta_1(D,B^*)}\Big\}\Big\},
\]
where the two inner $\max$ are taken over all $l=1,2,\dots,k-1$,
and over all $(l+1)\times l$ submatrix $D$ of $B$, for the first
$\max$, and of $B^*$, for the second.

It follows from our main theorem that there is a $k\times k$
matrix $B$ with $\pm 1$ entries for which
\begin{equation}\label{eq:formula}
\sqrt{2k}\max\Big\{\max\Big\{\frac{\Delta_2(D,B)}{\Delta_1(D,B)}\Big\},\max\Big\{\frac{\Delta_2(D,B^*)}{\Delta_1(D,B^*)}\Big\}\Big\}
\end{equation}
is bounded by a constant independent of $k$ and for such a matrix
this gives the splitting constant.

This of course gives an algorithm (still not very efficient) for
searching for a Kashin splitting, but more importantly, it suggest
that there might be an algebraic or combinatorial method of
finding an explicit splitting. The formula (\ref{eq:formula})
gives an explicit criterion for deciding whether a $\pm 1$ matrix
produce such a splitting.

\section{The main result}\label{main}

We shall denote by $\|a\|_p$ the $\ell_p^k$ norm of
$a=(a_1,a_2,\dots,a_n)$, $\|a\|_p=\Big(\sum_{i=1}^n
|a_i|^p\Big)^{1/p}$ (notice the difference with $\|x\|_{L_p^n}$
defined earlier). Denote by $S^{k-1}$ the Euclidean unit sphere in
$\R^k$ and by $B_p^k$ the unit ball of $\ell_p^k$, $0<
p\le\infty$. Given $a=(a_1,a_2,\dots,a_n)\in \R^k$ denote
\[
E_p(a)={\rm Ave}\Big( k^{-1}\sum_{j=1}^k|\sum_{i=1}^k a_i
\e_{i,j}|^p\Big)^{1/p}
\]
where the average is taken over all sequences of signs
$\{\e_{i,j}\}$. As is well known $2^{-1/2}\le E_1(a)\le E_2(a)\le
1$ (See \cite{sz} for the stated explicit lower bound, we only
need some absolute positive lower bound which follows from
Khinchine's inequality.) In the sequel $P$ denotes the natural
probability measure on $\{-1,1\}^{k^2}$ and the general element
in this probability space is denoted by $\{\e_{i,j}\}_{i,j=1}^k$.
We begin with two concentration inequalities.

\begin{lm}\label{lm:concentration} There is an absolute positive constant
$\eta$ such that for all $k$, all $a\in\R^n$ and all $0<C<\infty$,
\begin{equation}\label{eq:conc1}
P\Big(|k^{-1}\sum_{j=1}^k|\sum_{i=1}^ka_i\e_{i,j}|-E_1(a)|>CE_1(a)\Big)\le
e^{-\eta C^2 k}
\end{equation}
and
\begin{equation}\label{eq:conc2}
P\Big(|(k^{-1}\sum_{j=1}^k(\sum_{i=1}^ka_i\e_{i,j})^2)^{1/2}-E_2(a)|>C\Big)\le
e^{-\eta C^2 k}.
\end{equation}
\end{lm}

\pf Let $f,g:\R^{k^2}\to \R$ be the functions defined by
\[
f(x)=k^{-1}\sum_{j=1}^k|\sum_{i=1}^ka_i x_{i,j}| \ \ {\rm{ and }}\
\ g(x)= (k^{-1}\sum_{j=1}^k(\sum_{i=1}^ka_i x_{i,j})^2)^{1/2}
\]
for $x=\{x_{i,j}\}_{i,j=1}^k$. Both functions are convex and
Lipschitz with constant $\|a\|_2$ with respect to the
$\ell_2^{k^2}$ norm. The later statement can be proved by
computing the norm of the gradients of the functions at their
points of differentiability. The analogue inequality to
(\ref{eq:conc1}) where $E_1(a)$ is replaced with the median of
$f$ follows from the main result of \cite{ta}. That the median
can be replaced with the mean is simple, well known and can be
found e.g. in \cite{ms86} Proposition V.4. Inequality
(\ref{eq:conc2}) is dealt with similarly.
\endpf

For $1\le l\le k$ denote
\[
F_l^k=\{a=(a_1,a_2,\dots,a_k)\in B_2^k\
;\ a_i\not=0 \ {\rm{ for \ at \ most}}\ l \ {\rm{ values\ of}}\ i
\}.
\]
We now extend the concentration inequalities of Lemma
\ref{lm:concentration} to  simultaneous inequalities for the sets
$F_l^k$.

\begin{pr}\label{pr:smallsupp} There are absolute positive constants $\eta$ and
$\alpha$ such that for all $k$ and $l\le \alpha k$,
\[
P\Big(|k^{-1}\sum_{j=1}^k|\sum_{i=1}^ka_i\e_{i,j}|-E_1(a)|>E_1(a)/4,
\ {\rm{ for \ some\ }}a\in F_l^k\Big)\le e^{-\eta k}
\]
and
\[
P\Big(|(k^{-1}\sum_{j=1}^k(\sum_{i=1}^ka_i\e_{i,j})^2)^{1/2}-E_2(a)|>E_2(a)/4,
\ {\rm{ for \ some\ }}a\in F_l^k\Big)\le e^{-\eta k}.
\]
\end{pr}

\pf Given two sets $A, B$ in a linear space we denote by $N(A,B)$
the minimal number of shifts of $B$ whose union cover $A$.

For $\sigma\subseteq\{1,2,\dots,k\}$ denote
$F_{\sigma}^k=\{(a_1,a_2,\dots,a_k)\in B_2^k\ ;\ a_i=0 \ {\rm
{for}} \ i\notin \sigma\}$. Fix some $1\le l\le k$ then by the
usual volume estimates (see e.g. \cite{ms86}) for each subset
$\sigma\subseteq\{1,2,\dots,k\}$ of cardinality $l$ and for all
$0<\delta<1$, $N(F_{\sigma}^k,\delta B_2^k)\le (2\delta^{-1})^l$.
It follows form Lemma \ref{lm:concentration} that for some
absolute $\eta>0$
\[
P\Big(|k^{-1}\sum_{j=1}^k|\sum_{i=1}^ka_i\e_{i,j}|-E_1(a)|>E_1(a)/8,
\ {\rm{ for \ some\ }}a\in \calN \Big)\le e^{l\log(2/\delta)-\eta
k}
\]
where $\calN$ is some $\delta$-net in  $F_{\sigma}^k$. If $\delta$
is a small enough positive universal constant and $l/k$ is small
enough with respect to the universal constants $\delta$ and
$\eta$ (so that $e^{l\log(2/\delta)-\eta k}<e^{-\eta k/2}$, it
now follows by successive approximation (see e.g. \cite{ms86})
that
\begin{equation}\label{eq:b_sigma}
P\Big(|k^{-1}\sum_{j=1}^k|\sum_{i=1}^ka_i\e_{i,j}|-E_1(a)|>E_1(a)/4,
\ {\rm{ for \ some\ }}a\in F_{\sigma}^k \Big)\le e^{-\eta k/2}
\end{equation}

Put $\alpha=l/k$, assume $\alpha\le 1/2$ and also small enough for
(\ref{eq:b_sigma}) to hold. Notice that, by Stirling formula, the
number of subsets of $\{1,2,\dots,k\}$ of cardinality $l$ can be
evaluated as
\[
{k\choose l}\le e^{3k\alpha\log \frac{1}{\alpha}}.
\]

It follows from (\ref{eq:b_sigma}) that
\[
P\Big(|k^{-1}\sum_{j=1}^k|\sum_{i=1}^ka_i\e_{i,j}|-E_1(a)|>E_1(a)/4,
\ {\rm{ for \ some\ }}a\in F_l^k\Big)\le
e^{3k\a\log\frac{1}{\a}-\eta k/2}.
\]
Finally, if $\a\log\frac{1}{\a}<\eta/12$ the last quantity is
less than $e^{-\eta k/4}$ which finishes the proof of the first
assertion. The second is proved very similarly.
\endpf

\begin{lm}\label{lm:1} Let $a=(a_1,a_2,\dots,a_k)$ be a norm one
vector in $\ell_2^k$ and let $0<\gamma<1$. Assume
$k^{-1/2}\sum_{i=1}^k|a_i|\le\gamma$ then
\[
(\sum_{i=l}^k (a_i^*)^2)^{1/2}\le l^{-1}\gamma\sqrt{k(k-l+1)}
\]
for all $1\le l\le k$, where $\{a_i^*\}$ denotes the decreasing
rearrangement of $\{|a_i|\}$.
\end{lm}

\pf For each $1\le l\le k$,  \ $(k-l+1)(a_l^*)^2\ge \sum_{i=l}^k
(a_i^*)^2$. It follows that
\[
\gamma\ge k^{-1/2}\sum_{i=1}^k|a_i|\ge lk^{-1/2}a_l^*\ge
l(k(k-l+1))^{-1/2}(\sum_{i=l}^k (a_i^*)^2)^{1/2}
\]
from which the conclusion follows.
\endpf

For $0<\gamma<1$, $k\in\N$, denote
\[
A_{\gamma}^k=\{a\in S^{k-1}\ ;\
k^{-1/2}\sum_{i=1}^k|a_i|\le\gamma\}=\Big\{a\in S^{k-1}\ ;\
\frac{\|a\|_{L_1^k}}{\|a\|_{L_2^k}}\le \gamma\Big\}.
\]
Next we extend the concentration inequalities to the sets
$A_{\gamma}^k$.

\begin{pr}\label{pr:main}There are  absolute constants
$0<\gamma<1$ and $\eta>0$ such that for all $k$,
\[
P\Big(|k^{-1}\sum_{j=1}^k|\sum_{i=1}^ka_i\e_{i,j}|-E_1(a)|>E_1(a)/2,
\ {\rm{ for \ some\ }}a\in A_{\gamma}^k\Big)\le e^{-\eta k}
\]
and
\[
P\Big(|(k^{-1}\sum_{j=1}^k(\sum_{i=1}^ka_i\e_{i,j})^2)^{1/2}-1|>1/2,
\ {\rm{ for \ some\ }}a\in A_{\gamma}^k\Big)\le e^{-\eta k}.
\]
\end{pr}

\pf Notice first that by the usual $\e$-net considerations
starting with Lemma \ref{lm:concentration}, there are some
absolute $C<\infty$ and $\eta>0$ such that
\[
P\Big(k^{-1}\sum_{j=1}^k|\sum_{i=1}^ka_i\e_{i,j}|>C\|a\|_2, \
{\rm{ for \ some\ }}a\in \R^k\Big)\le e^{-\eta k}
\]
and
\[
P\Big((k^{-1}\sum_{j=1}^k(\sum_{i=1}^ka_i\e_{i,j})^2)^{1/2}>C\|a\|_2,
\ {\rm{ for \ some\ }}a\in \R^k\Big)\le e^{-\eta k}.
\]
(Actually, the first assertion follows trivially from the second.)
Let $l=[\alpha k]$ where $\alpha$ is the constant from Proposition
\ref{pr:smallsupp}. Now choose $\gamma>0$ such that, putting
$\delta=l^{-1}\gamma\sqrt{k(k-l+1)}$, \
$3C\delta<\frac{1}{4\sqrt{2}}$.

By Lemma \ref{lm:1}, each $a\in A_{\gamma}^k$ can be split as
$a=b+c$ with $b\in F_l^k$ and $\|c\|_2\le \delta$. Let
$\{\e_{i,j}\}$ be such that both
\begin{equation}\label{eq:dev}
|k^{-1}\sum_{j=1}^k|\sum_{i=1}^kb_i\e_{i,j}|-E_1(b)|\le E_1(b)/4
\end{equation}
and
\begin{equation}\label{eq:bound}
k^{-1}\sum_{j=1}^k|\sum_{i=1}^kc_i\e_{i,j}|\le C\|c\|_2.
\end{equation}
Then the
condition on $\delta$ implies that
\[
|k^{-1}\sum_{j=1}^k|\sum_{i=1}^ka_i\e_{i,j}|-E_1(a)|\le E_1(a)/2.
\]
Since, by Proposition \ref{pr:smallsupp} and the first paragraph
of this proof,  the probability that at least one of inequalities
(\ref{eq:dev}) or (\ref{eq:bound}) does not hold is less than
$e^{-\eta k}$, we get the first assertion of the Proposition
(with a different absolute $\eta$). The second assertion is
proved very similarly.
\endpf

Given signs $\{\e_{i,j}\}_{i,j=1}^k$ we shall denote by $B$ the
$k\times k$ matrix with entries $\{\e_{i,j}\}_{i,j=1}^k$ and by
$A=[\sqrt{k}I,B]$ the $k\times 2k$ matrix whose first $k$ columns
form $\sqrt{k}I_k$ and the last $k$ columns form $B$. We shall
also denote $\bar A=[-B^*,\sqrt{k}I]$ with the obvious meaning
(where $B^*$ is the transpose of $B$). Note that the rows span of
$A$ and of $\bar A$ form orthogonal subspaces of $L_2^{2k}$. We
are now ready to state and prove our main result.

\begin{thm}\label{thm:main}
For some absolute $\eta>0$ and $C<\infty$ and for all $k$ there
are signs $\{\e_{i,j}\}_{i,j=1}^k$ such that for all $a\in
S^{k-1}$
\[
C^{-1}\le \|aA\|_{L_1^{2k}}\le \|aA\|_{L_2^{2k}}\le C
\]
and
\[
C^{-1}\le \|a\bar A\|_{L_1^{2k}}\le \|a\bar A\|_{L_2^{2k}}\le C.
\]
Moreover, this holds with probability larger than $1-e^{-\eta k}$.
\end{thm}

\pf As in the beginning of the proof of Proposition
\ref{pr:main}, it follows from Lemma \ref{lm:concentration} that
for some absolute $C$ and $\eta$ and with probability at least $
1-e^{-\eta k}$,
\[
k^{-1}\sum_{j=1}^k|\sum_{i=1}^ka_i\e_{i,j}|\le
(k^{-1}\sum_{j=1}^k(\sum_{i=1}^ka_i\e_{i,j})^2)^{1/2}\le C
\]
for all $a\in S^{k-1}$. Of course the same holds also if we
replace $\e_{i,j}$ with $-\e_{j,i}$ everywhere. It follows easily
that, with probability $\> 1-e^{-\eta k}$,
\[
\|aA\|_{L_1^{2k}}\le \|aA\|_{L_2^{2k}}\le
\frac{(1+C^2)^{1/2}}{\sqrt{2}}
\]
and
\[
\|a\bar A\|_{L_1^{2k}}\le \|a\bar A\|_{L_2^{2k}}\le
\frac{(1+C^2)^{1/2}}{\sqrt{2}}
\]
 for all $a\in S^{k-1}$.

For the lower bound let $\gamma$ be the constant from Proposition
\ref{pr:main}. Then, with probability $\> 1-e^{-\eta k}$,
\[
k^{-1}\sum_{j=1}^k|\sum_{i=1}^ka_i\e_{i,j}|>\frac{1}{2\sqrt{2}}
\]
for all $a\in A_{\gamma}^k$ (and the same holds with $-\e_{j,i}$
instead of $\e_{i,j}$). For the other $a\in S^{k-1}$
\[
k^{-1/2}\sum_{i=1}^k|a_i|>\gamma.
\]
It follows easily that, with probability at least $1-e^{-\eta k}$,
\[
\|aA\|_{L_2^{2k}}\ge \|aA\|_{L_1^{2k}}\ge
\min\{\frac{\gamma}{2},\frac{1}{4\sqrt{2}}\}
\]
and
\[ \|a\bar A\|_{L_2^{2k}}\ge \|a\bar A\|_{L_1^{2k}}\ge
\min\{\frac{\gamma}{2},\frac{1}{4\sqrt{2}}\}
\]
for all $a\in S^{k-1}$.
\endpf

Recall the definition of $\Delta_p(D,B)$ appearing in the
Introduction.

\begin{co}\label{co:determinants}
There is a constant $C<\infty$ such that for all $k$ there is a
$k\times k$ matrix $B$ with $\pm 1$ entries such that
\begin{equation}\label{eq:formula1}
\sqrt{2k}\max\Big\{\max\Big\{\frac{\Delta_2(D,B)}{\Delta_1(D,B)}\Big\},
\max\Big\{\frac{\Delta_2(D,B^*)}{\Delta_1(D,B^*)}\Big\}\Big\}\le C
\end{equation}
where the first inner $\max$ is taken over all $l=1,2,\dots,k-1$,
and over all $(l+1)\times l$ submatrix $D$ of $B$ for which the
denominator is not zero, while the second inner $\max$ is taken
over all $l=1,2,\dots,k-1$, and over all $(l+1)\times l$
submatrix $D$ of $B^*$. \hfill\break Moreover, the left hand side
of (\ref{eq:formula1}) is equal to
$\max\{C_{2k}(E),C_{2k}(E^{\bot}\}$, where $E$ is the span of the
rows of $[\sqrt{k}I,B]$.
\end{co}

\pf By Theorem \ref{thm:main} we only need to address the
``Moreover" part. This follows easily from section 2.5 in
Anderson's \cite{a}, in which it is shown, in our notations, that
for a $k\times 2k$ matrix $A$ and for $E$ being the span of its
rows,
\[
C_{2k}(E)=\sqrt{\frac{2k}{k+1}}\max\frac{\Big(\frac{1}{k+1}\sum_{j\notin
C}|\det D^{+j}|^2\Big)^{1/2}}{\frac{1}{k+1}\sum_{j\notin C}|\det
D^{+j}|}
\]
where the max is taken over all $k\times(k-1)$ submatrices
$D=D_{\{1,\dots,k\},C}$ for which the denominator does not vanish.

\endpr

As we remarked in the introduction Corollary \ref{co:determinants}
gives an explicit criterion for deciding whether a $\pm 1$ matrix
gives a good Kashin splitting and the asurance that there are
(many) such matrices. We hope this will help in a search for an
explicit construction of a Kashin splitting.

\medskip

\noindent{\bf Remark:} \label{co:intersection} It is easy to see
that Theorem \ref{thm:main} implies that for each positive integer
$k$, with probability larger than $1-e^{-\eta k}$, a $k\times k$
matrix $B$ with independent $\pm 1$ entries satisfies the
following: Letting $K_1$ be $\{x\ ;\ Bx\in \sqrt{k}B_1^k\}$ and
$K_2$ be the symmetric convex hull of $\sqrt{k}$ times the
canonical unit vector basis in $\R^k$ ($=\sqrt{k}B_1^k$), then
$K_1\cap K_2$ lies between two universal multiples of the
Euclidean unit ball, $B_2^k$.

\section{Some related results and remarks}

Kashin also proved that for {\em any} $0<\lambda <1$ and any $n$
there is a $[\lambda n]$-dimensional subspace $E$ of $L_1^n$
which is $C(\lambda)$-isomorphic to a Hilbert space where
$C(\lambda)$ depends only on $\lambda$ (actually, he proved that
$C_n(E)\le C(\lambda)$). A similar statement holds for almost
isometries although we need to replace ``any $0<\lambda <1$" with
``some $0<\lambda <1$": For every $\e>0$ there is a
$\lambda=\lambda(\e)$ such that for any $n$ there is a subspace of
$L_1^n$ of dimension at least $\lambda n$ which is
$(1+\e)$-isomorphic to a Hilbert space. This is proved in
\cite{flm} (and, without stating it explicitly,  already in
\cite{mi}). The proofs are again probabilistic and we are far from
having any explicit embeddings.

What about embeddings given by span of rows of matrices whose
entries take values in some small set of values? It is not very
hard to see that a similar proof to the one here gives, for any
$0<\lambda <1$ and any $n$, a subspace $E$ of $L_1^n$ on which
$C_n(E)\le C(\lambda)$ and which is spanned by vectors whose
entries are taken from a four point set (actually, the set
consists of $0, \pm 1$ and one other specific value, only the
$\pm 1$ are chosen randomly). Moreover, there is a corresponding
determinantal formula for determining whether a space from this
collection satisfies $C_n(E)\le C(\lambda)$. These subjects will
be detailed in a forthcoming MSc thesis of Boris Levant written
at the Weizmann Institute.

It is also possible to find, for every $0<\lambda<1$, a good
Hilbertian subspace of $L_1^n$ of dimension $[\lambda n]$ spanned
by rows of a $[\lambda n]\times n$ matrix with $\pm 1$ entries.
This follows from the main result of \cite{sch} where a similar
statement with {\em some} $0<\lambda<1$ instead of {\em every}
$0<\lambda<1$ is proved, together with the method of \cite{js03}
where it is shown that whenever $E$ is a $k$-dimensional subspace
of $L_1^n$ then, for all $a>1$, the restriction operator onto some
$[ak]$ of the $n$ coordinates is a $C$-isomorphism when
restricted to $E$ and where $C$ depends on $k/n$ only. The proof
in \cite{sch} uses a concentration inequality similar to the one
in the first part of Lemma \ref{lm:concentration}. Using a
variation on the second part of that lemma as well (and the
restriction method of \cite{js03}) one can get a bit more.

\begin{pr}\label{pr:restriction}
For all $0<\lambda<1$ and all $n$ there is a $[\lambda n]\times n$
matrix $A$ with $\pm 1$ entries such that for all $a\in
S^{[\lambda n]-1}$
\[
C^{-1}(\lambda)\le \|aA\|_{L_1^{n}}\le \|aA\|_{L_2^{n}}\le
C(\lambda)
\]
Where $C(\lambda)$ depends on $\lambda$ only.
\end{pr}
The details of the proof will be given in Levant's thesis.

\bigskip

\noindent {\bf Remark:} Going back to the search for an explicit
$\pm 1 $ matrix $B$ for which the span of the rows of
$[\sqrt{k}I,B]$ gives a good Kashin splitting, a first candidate
to look for is the Walsh matrix. However, it is easy to see that
this is not the case. Assume $k=2^t$; reindex the columns
$1,\dots,k$ as $\{-1,1\}^t$ and the rows by $
\{\sigma\}_{\sigma\subseteq\{1,\dots,t\}}$ and let the
$\sigma,\e$ term  of the matrix $B$ be
$W_\sigma(\e)=\prod_{i\in\sigma}\e_i$. Consider the vector of
coefficients $a=(a_\sigma)$ where $a_\sigma$ is $1$ whenever
$\sigma$ is a subset of $\{1,\dots,t/2\}$ (assuming $t$ is even)
and $0$ otherwise. Then it is not hard to see that
$\|\sqrt{k}a\|_{L_2^k}=\|aB\|_{L_2^k}=k^{1/4}$, while
$\|\sqrt{k}a\|_{L_1^k}=\|aB\|_{L_1^k}=1$.

\bigskip

The method of the proof of the main theorem may be useful for
other applications. The idea of the proof was that we split the
sphere $S^{k-1}$ into two sets. On one of them the $L_1^k$ and
$L_2^k$ are well equivalent and the other one ($A_\gamma^k$) is
``small". Of course the measure of $A_\gamma^k$ is basically
known and is very small. This estimate was not good enough for
our purposes and we needed another measure of ``smallness" (which
was $A_\gamma^k\subset F_l^k+\delta B_2^k$ for some (not too
small) $l$ and (small) $\delta$). There is another measure of
smallness that follows easily from the proof here and may be
useful elsewhere. Again, it was not good enough for our purposes.
Recall that $A_{\gamma}^k=\{(a_1,a_2,\dots,a_k)\in S^{k-1}\ ;\
k^{-1/2}\sum_{i=1}^k|a_i|\le\gamma\}$.

\begin{pr}\label{pr:nets}
Let $0<\gamma<1$ and $k\in\N$. Then, for all $\e>4\gamma$,
\[
N(A_{\gamma}^k,\e B_2^k)\le e^{\frac{6\gamma
k}{\e}(\log\frac{\e}{2\gamma}+\log\frac{4}{\e})}.
\]
\end{pr}

\pf Since for any $\sigma\subseteq\{1,2,\dots,k\}$ of cardinality
$l$ and for all $0<\delta<1$, $N(F_{\sigma}^k,\delta B_2^k)\le
(2\delta^{-1})^l$. and since the number of subsets of
$\{1,2,\dots,k\}$ of cardinality $l\le k/2$ can be evaluated as
\[
{k\choose l}\le e^{3l\log \frac{k}{l}},
\]
it follows that
\begin{equation}\label{eq:blk}
N(F_l^k,\delta B_2^k)\le e^{3l(\log \frac{k}{l}+\log
\frac{2}{\delta})}.
\end{equation}

By Lemma \ref{lm:1}, $A_\gamma^k\subset F_l^k+\delta B_2^k$ with
$\delta=l^{-1}\gamma\sqrt{k(k-l+1)}$. It follows that
\[
N(A_\gamma^k,2\delta B_2^k)\le e^{3l(\log \frac{k}{l}+\log
\frac{2}{\delta})}.
\]
Letting $\delta=\e/2$ and $l=[2\gamma\sqrt{k(k-l+1)}/\e]\le
2\gamma k/\e$, we get for $\gamma<\e/4$ (to ensure $l<k/2$),
\[
N(A_\gamma^k,2\delta B_2^k)\le e^{\frac{6\gamma
k}{\e}(\log\frac{\e}{2\gamma}+\log\frac{4}{\e})}.
\]
\endpf

\bigskip

\bigskip

\noindent Gideon Schechtman\newline
Department of Mathematics\newline
Weizmann Institute of Science\newline
Rehovot, Israel\newline
E-mail: gideon@wisdom.weizmann.ac.il


\begin{thebibliography}{99}



\bibitem[An]{a}
Anderson, G. W., Integral Ka\v{s}in splittings, to appear in the
Israel J. of Math.




\bibitem[FLM]{flm}
Figiel, T.; Lindenstrauss, J.; Milman, V. D., The dimension of
almost spherical sections of convex bodies. Acta Math. 139
(1977), no. 1-2, 53--94.


\bibitem[Ka]{ka}
Kashin, B., Section of some finite-dimensional sets and classes
of smooth functions. Izv. Acad. Nauk. SSSR 41 (1977), 334--351.
(Russian).






\bibitem[JS]{js03} Johnson, W. B.; Schechtman, Very tight
embeddings of subspaces of $L_p$, $1\le p<2$, into $\ell_p^n$,
Geom. Funct. Anal., to appear

\bibitem[Mi]{mi}
Milman, V. D., A new proof of A. Dvoretzky's theorem on
cross-sections of convex bodies. (Russian) Funkcional. Anal. i
Prilo\v zen. 5 (1971), no. 4, 28--37.

\bibitem[MS]{ms86} Milman,
V. D. and Schechtman, G., Asymptotic theory of finite-dimensional
normed spaces, Lecture Notes in Mathematics, 1200,
Springer-Verlag, Berlin, 1986.


\bibitem[Sc]{sch}
Schechtman, G., Random embeddings of Euclidean spaces in sequence
spaces. Israel J. Math. 40 (1981), no. 2, 187--192.

\bibitem[Sz]{sz}
Szarek, S. J., On the best constants in the Khinchin inequality.
Studia Math. 58 (1976), no. 2, 197--208.

\bibitem[Ta]{ta}
Talagrand, M., An isoperimetric theorem on the cube and the
Kintchine-Kahane inequalities. Proc. Amer. Math. Soc. 104 (1988),
no. 3, 905--909.



\end{thebibliography}
\end{document}